\documentclass[article]{jotart}
\usepackage{amsmath,amstext, empheq, extpfeil,  
amsbsy,marvosym,fancyhdr,graphics,amscd,amsfonts,latexsym,delarray,stackrel,lineno,color,cite,appendix,xcolor,url,hyperref}
\usepackage{amssymb}
\usepackage{wasysym}
\usepackage{subfig}
\usepackage{srcltx}
\usepackage{mathrsfs}
\usepackage{float} 
\usepackage[all]{xy}
\usepackage{t1enc}
\usepackage{mathrsfs}
\usepackage{pifont}
\usepackage{mathpple}
\usepackage[T1]{fontenc}

\definecolor{Green}{rgb}{0,1,0}
\definecolor{Blue}{RGB}{0,0,191}
\definecolor{mathmodecolor}{RGB}{0,102,0}
\definecolor{keywordcolor}{RGB}{0,51,151}
\definecolor{sourcebackgroundcolor}{RGB}{255,247,223}
\definecolor{unixagred}{RGB}{255,0,0}
\definecolor{lightgray}{RGB}{191,191,191}
\definecolor{green}{RGB}{1,191,191}

\newcommand*\patchAmsMathEnvironmentForLineno[1]{%
  \expandafter\let\csname old#1\expandafter\endcsname\csname #1\endcsname
  \expandafter\let\csname oldend#1\expandafter\endcsname\csname end#1\endcsname
  \renewenvironment{#1}%
     {\linenomath\csname old#1\endcsname}%
     {\csname oldend#1\endcsname\endlinenomath}}%
\newcommand*\patchBothAmsMathEnvironmentsForLineno[1]{%
  \patchAmsMathEnvironmentForLineno{#1}%
  \patchAmsMathEnvironmentForLineno{#1*}}%
\AtBeginDocument{%
\patchBothAmsMathEnvironmentsForLineno{equation}%
\patchBothAmsMathEnvironmentsForLineno{align}%
\patchBothAmsMathEnvironmentsForLineno{flalign}%
\patchBothAmsMathEnvironmentsForLineno{alignat}%
\patchBothAmsMathEnvironmentsForLineno{gather}%
\patchBothAmsMathEnvironmentsForLineno{multline}%
}

\theoremstyle{proclaim}
\newtheorem{theorem}{Theorem}[section]
\newtheorem{lemma}[theorem]{Lemma}
\newtheorem{corollary}[theorem]{Corollary}

\newtheorem{guess}[theorem]{Conjecture}
\theoremstyle{statement}
\newtheorem{remark}[theorem]{Remark}
\newtheorem{definition}[theorem]{Definition}

\newtheorem{fact}[theorem]{Fact}
\theoremstyle{fancyproclaim}

\numberwithin{equation}{section}

\def\qd{^-\!\!\!\!\!d}

\def\F{{\mathbb F}}

\def\qqq{\,,\quad \forall}

\def\bm2{{\rm \B mod^2}}
\def\b2{{\rm \B mod^{\mathfrak s}}}

\def\GL{{\rm GL}}

\def\Mod{{\rm Mod}}

\def\Spec{{\rm Spec}}

\def\Tr{{\rm Tr}}

\def\A{{\mathbb A}}
\def\C{{\mathbb C}}
\def\F{{\mathbb F}}

\def\Q{{\mathbb Q}}
\def\R{{\mathbb R}}

\def\Z{{\mathbb Z}}
\def\H{{\mathbb H}}

\def\B{{\mathbb B}}

\def\fourier{\F}

\def\Tr{{\rm Tr}}

\def\cA{{\mathcal A}}

\def\cD{{\mathcal D}}

\def\cH{{\mathcal H}}

\def\cS{{\mathcal S}}

\def\part{\partial}

\def\mapE{{\mathfrak E}}

\newcommand{\ie}{{\it i.e.\/}\ }

\def\sin{{{\rm sin}}}

\def\H{{\mathbb H}}

\def\urep{\vartheta_{\rm a}}
\def\vrep{\vartheta_{\rm m}}

\definecolor{trust}{rgb}{0,1,1}

\def\noncommutative geometry{{noncommutative geometry }}

\def\qqq{\,,\quad \forall}

\begin{document}
\issueinfo{00}{0}{0000} 
\commby{Editor}
\pagespan{101}{110}
\date{Month dd, yyyy}
\revision{Month dd, yyyy}
\title[The scaling Hamiltonian]{The scaling Hamiltonian}
\dedicatory{Dedicated to Dan Voiculescu with admiration}
\author[Alain Connes {\protect \and} Caterina Consani]{Alain Connes {\protect \and} Caterina Consani}
\address{Alain Connes, Coll\`ege de France, IHES and Ohio State University}
\email{alain@connes.org}
\address{Caterina Consani, Department of Mathematics, The Johns Hopkins University, Baltimore MD, 21218, USA}
\email{kc@math.jhu.edu}
\begin{abstract} 
We first explain the link between the  Berry-Keating Hamiltonian  and the spectral realization of zeros of the Riemann zeta function $\zeta$ of \cite{Co-zeta}, and why there is no conflict at the semi-classical level between the ``absorption"  picture of  \cite{Co-zeta} and the semiclassical ``emission" computations of \cite{BK, BKe}, while the minus sign manifests itself in the Maslov phases. We then use the quantized calculus to analyse the recent attempt of X.-J.~Li at proving  Weil's positivity, and understand its limit. We then propose an operator theoretic semi-local framework directly related to the Riemann Hypothesis.  
\end{abstract}
\maketitle

\section*{INTRODUCTION}

This special volume of the Journal of Operator Theory is a wonderful occasion to acknowledge the scientific contributions of Dan Voiculescu  both as a problem solver and a theory builder. 

The topic of the present paper can be described as a nagging thought on the power of Operator Theory in connection with a main unsolved problem in mathematics : the Riemann Hypothesis (RH). After the initial attempt of the first author \cite{Co-zeta}, that stemmed from a coincidental encounter with RH through quantum statistical mechanics, both authors (of the present paper) abandoned the hope of dealing with the problem in a purely analytic manner and shifted their attention, in search of another powerful tool : Algebraic Geometry. We refer to \cite{Crh,CCsurvey} for a  survey  of the geometric  space that our approach has unveiled (the Scaling Site) and to \cite{schemeF1} for the development of an ``absolute algebraic geometry", motivated by the discovery of this new site, and the need to formulate a Riemann-Roch formula on its square.

The nagging thought which is the subject of the present paper is that, after all, one might have given up too easily on analysis and the power of Hilbert space operators to obtain positivity statements. Indeed in  \cite{Co-zeta}, while the construction of a global Hilbert space spectral realization remained artificial by the use of Sobolev spaces, the construction of the semi-local Hilbert space operators is completely canonical and provides the subtle finite parts in the Riemann-Weil explicit formulas.   We first explain in Section \ref{sectscalhamil} the relation between the so-called Berry-Keating Hamiltonian\footnote{Their work was motivated by  extensive explanations of A. Connes at  IHES, as M. Berry and J. Keating  acknowledge in \cite{BK}} and the  original spectral realization of  \cite{Co-zeta}. The key points  are first why passing from an absorption spectrum (as in \cite{Co-zeta}) to an emission spectrum (as suggested  in \cite{BK,BKe}) does {\em not} introduce a mismatch in the semiclassical approximations, and moreover why the cutoff suggested in  \cite{BK,BKe} is equivalent, at the semiclassical level to the removal of the contribution of the ``white light" from the quantum system. We also point out that there is no issue about quantizing the $H=PQ$ Hamiltonian, and this is what was already done in  \cite{Co-zeta}, and that the minus sign is not completely eliminated in the cutoff suggested in  \cite{BK,BKe} and appears in the Maslov phases.

In Section \ref{sectsemilocal} we recall the semi-local trace formula proved in \cite{Co-zeta}  as formulated in \cite{CMbook}. The main technical result, recalled in  Lemma \ref{techn1},  computes the trace of the scaling action on the  Hilbert space $L^2(X_{S})$ canonically associated to a finite set $S\ni \infty$ of places of $\Q$. The space $X_{S}$ is an approximation of the ad\` ele class space of the rationals, obtained by trading the full  ad\` ele ring $\A_\Q$ for the finite product of the local fields $\Q_v$ attached to the places $v\in S$, and effecting its quotient by the group 
\begin{equation}\label{GL1QS}
\Q_S^*=\GL_1(\Q_S)= \{ \pm p_1^{n_1} \cdots p_k^{n_k} \, :\,  p_j
\in S \setminus\{ \infty \} \,,\, n_j\in \Z\}.
\end{equation}

In Section \ref{sectsonin} we describe the idea of X.-J.~Li as an attempt to prove Weil's positivity. We use the quantized calculus to analyse the cutoff that he proposes, motivated by Sonine spaces. We prove a key Lemma \ref{semiloclem} showing that his approach would have worked out, had a certain function $\phi(z)$   belonged to the 
Hardy space $\H^\infty$ of the half-plane $\C^+_{\frac 12}:=\{z\in \C\mid \Re(z)>1/2\}$. The function $\phi(z)$ is meromorphic in $ \C$ and of modulus $1$ on the critical line
 so that, at first sight, it
looks like an "inner function", but unfortunately it fails to be so and in particular  
 it is not bounded. It should have been clear from the start though,  that such a semi-local approach could not work without using a global Poisson Formula that dictates the normalizations of the principal values.
By ignoring all primes above some size (a natural  strategy)  one has normalized the principal value
to be 0 for them, thus the local normalization of the principal values for the primes in $S$ should play a key
role, and one should not be able to neglect it completely (in \cite{Co-zeta} the $2 \log \Lambda$ term accommodates such local changes). 

In Section \ref{sectoptheory} we formulate a conjecture on the power of the semi-local framework as a tool to prove Weil's positivity. We explain why the global Poisson Formula makes sense in the semi-local cases and provides a  conceptual explanation for the functions involved in the semi-local trace formula. We finally discuss the  merits of the geometric strategy compared to the analytic one outlined here.

\section{THE SCALING HAMILTONIAN}\label{sectscalhamil}

In the setting of \cite{Co-zeta} (see also Chapter 2 of \cite{CMbook}) the spectral realization is as an absorption
spectrum, \ie it appears as dark lines in a background given by white light.

\begin{figure}[H]
\includegraphics[scale=0.5]{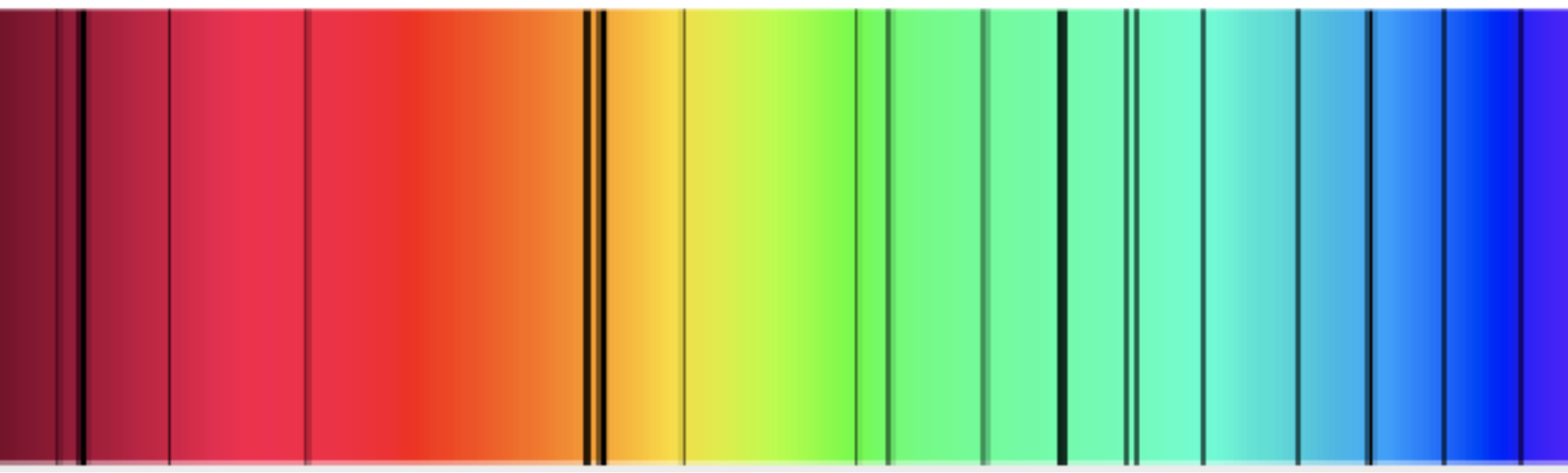}
\caption{Absorption Spectrum\label{figone}}
\end{figure} When disregarding the finite places to get the  simplest instance of the semi-local  picture, one obtains the action of the multiplicative group $\R_+^*$ by scaling on the Hilbert space $L^2(\R)$ of square integrable functions on the real line. This procedure corresponds in an obvious way to the ``quantization" of the Hamiltonian $H=PQ$. There is no mystery whatsoever in this quantization since in the unique irreducible representation of the Heisenberg commutation relations the operators $P$ and $Q$ correspond respectively (up to normalization) to $i\partial_x$ and multiplication by $x$, so that the self-adjoint operator $\frac 12(PQ+QP)$ corresponds to  the generator of the scaling flow which is $x\partial_x$ (to which one needs to add the constant $\frac 12$ in order to get the generator of a unitary group).   However, since the Hamiltonian is unbounded below, the associated quantum system given by the scaling action on
$L^2(\R)$ does not have discrete spectrum. In order to analyze it,   one performs a cutoff both in $q$  and $p$ spaces, at the same size $\Lambda$. Due to the normalization of the Fourier transform used in number theoretic contexts (\ie with phase factor $e^{-2\pi i pq}$), the Hamiltonian is $H=2\pi pq$ (see Chapter II, \S 3.2 of \cite{CMbook}). The portion of phase space in which the size of $H=2\pi pq$ is less than a given energy $E$ is shown in  Figure \ref{figone}.

\begin{figure}[H]
\begin{center}
\includegraphics[scale=0.8]{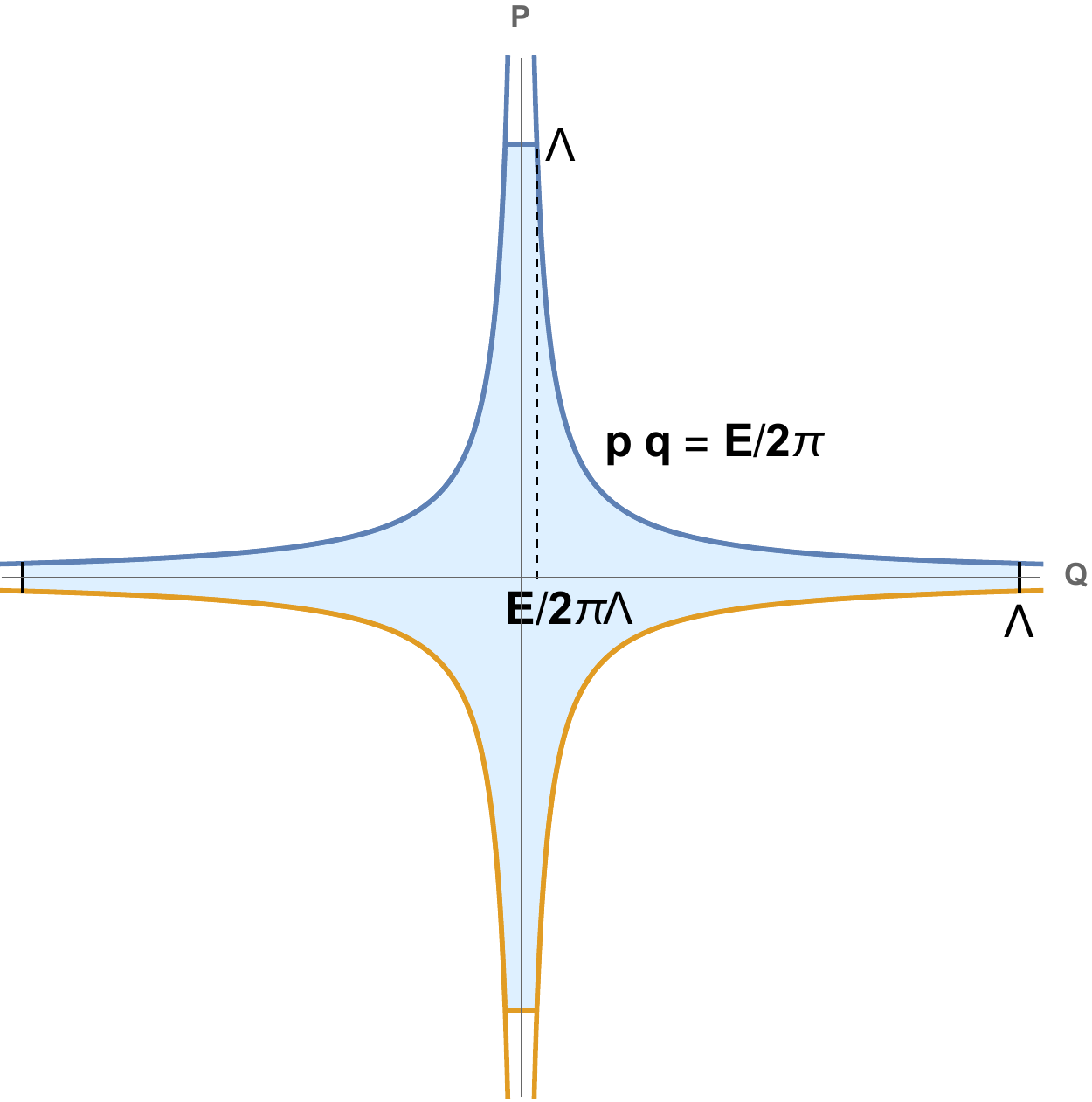}
\end{center}
\caption{\label{figone} }
\end{figure}

\begin{figure}[H]
\begin{center}
\includegraphics[scale=1]{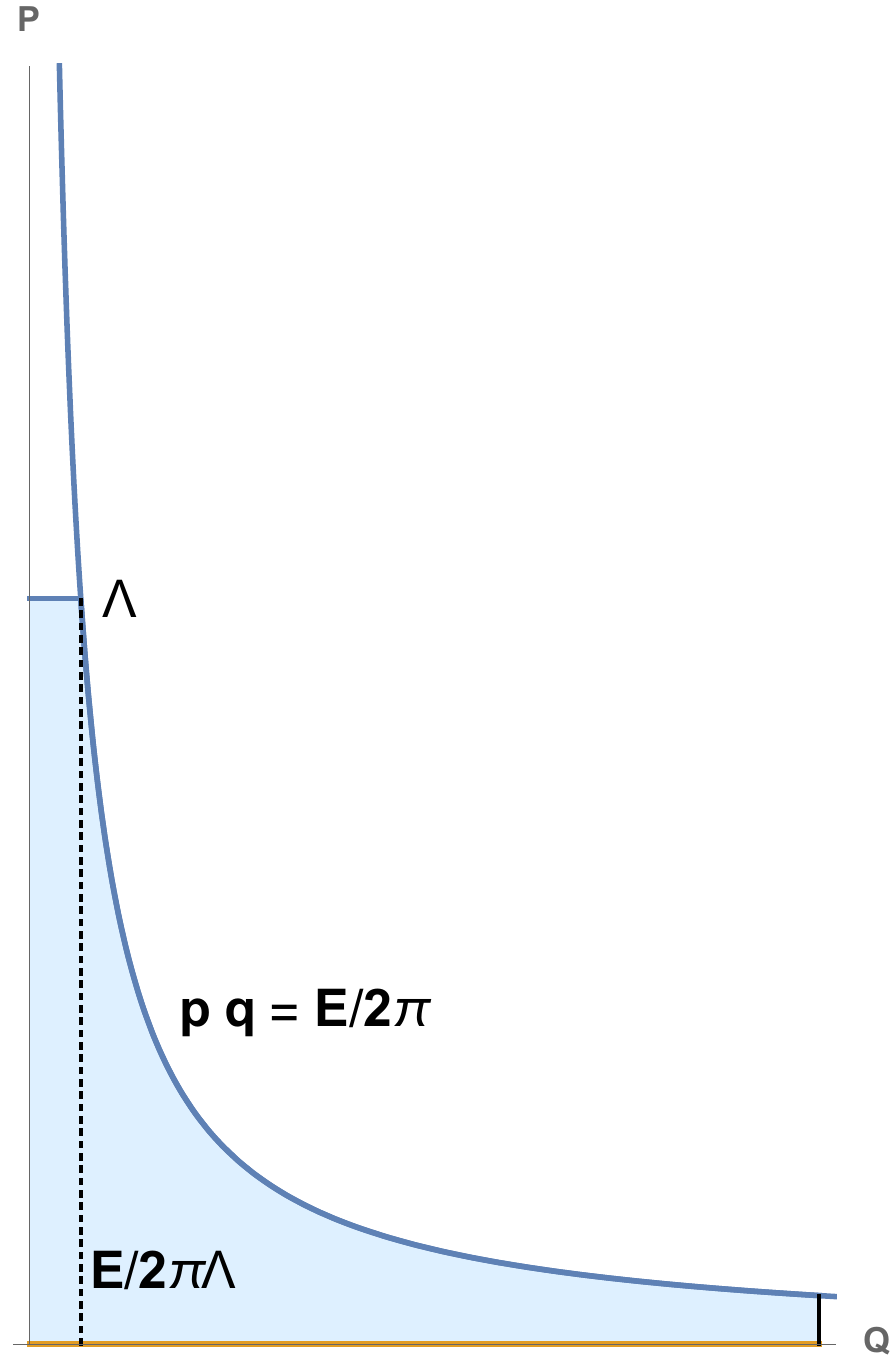}
\end{center}
\caption{\label{fig2} }
\end{figure}

\begin{figure}[H]
\begin{center}
\includegraphics[scale=0.6]{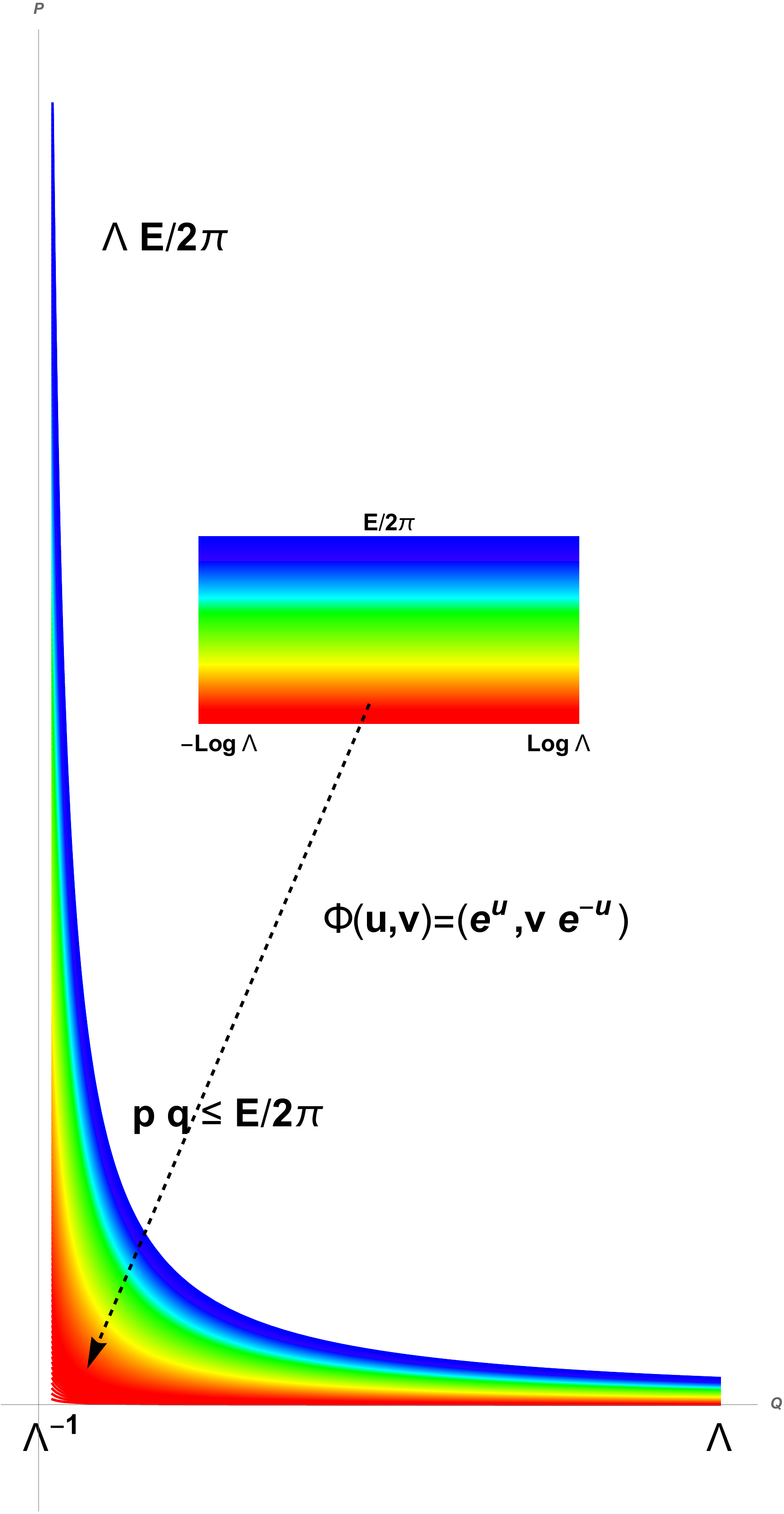}
\end{center}
\caption{\label{figthree} }
\end{figure}

As explained in details in  Chapter II, \S 3.2 of \cite{CMbook}, one restricts to  even functions to handle the Riemann zeta function, and 
 to the first quadrant in order to estimate the semiclassical count of the zeros of $\zeta$ with positive imaginary part less than $E$.  This procedure is depicted in  Figure \ref{fig2}.  In order to compare this semiclassical picture with white light,  one considers the model of white light given by the regular representation of the additive group $\R$. We adjust the cutoff so that  this corresponds, using the isomorphism  $u\mapsto \exp(u)$ of $\R$  with the multiplicative group $\R_+^*$, to 
the subgraph of the hyperbola $pq=\frac{E}{2\pi}$ restricted to the interval $[1/\Lambda,\Lambda]$ as shown in Figure \ref{figthree}.

The exponential $u\mapsto \exp(u)$ is obviously an isomorphism $\exp:\R\to \R_+^*$ and the following map is a symplectic isomorphism 
$$
\Phi(u,v):=(\exp(u),v \, \exp(-u)),
$$
since its Jacobian is given by 
$$
d\exp(u)\wedge d(v \, \exp(-u))=\exp(u) du \wedge (\exp(-u) dv-v \exp(-u) du)=du \wedge dv.
$$

 The symplectic isomorphism $\Phi$ transforms the phase space of the regular representation of $\R$ cutoff in the interval $[-\log\Lambda,\log\Lambda]$ and in  a range of energy (the dual variable to $\log(q)$) in the interval $[0,E/2\pi]$, into the portion of the phase space $(q,p)$ shown in Figure \ref{figthree}. We refer to
 \cite{Co-zeta} and Chapter II, \S 3.2 of \cite{CMbook} for the understanding of careful normalizations. 
There are two main differences between the cutoff phase space of  Figure \ref{fig2} and the white light of Figure \ref{figthree}. Firstly, there is a piece of Figure \ref{figthree} which is not covered by Figure \ref{fig2}, it extends from $p=\Lambda$ to $p=\Lambda E/2\pi$ and is represented in dark blue in Figure \ref{figfour} :
\begin{figure}[H]
\begin{center}
\includegraphics[scale=0.9]{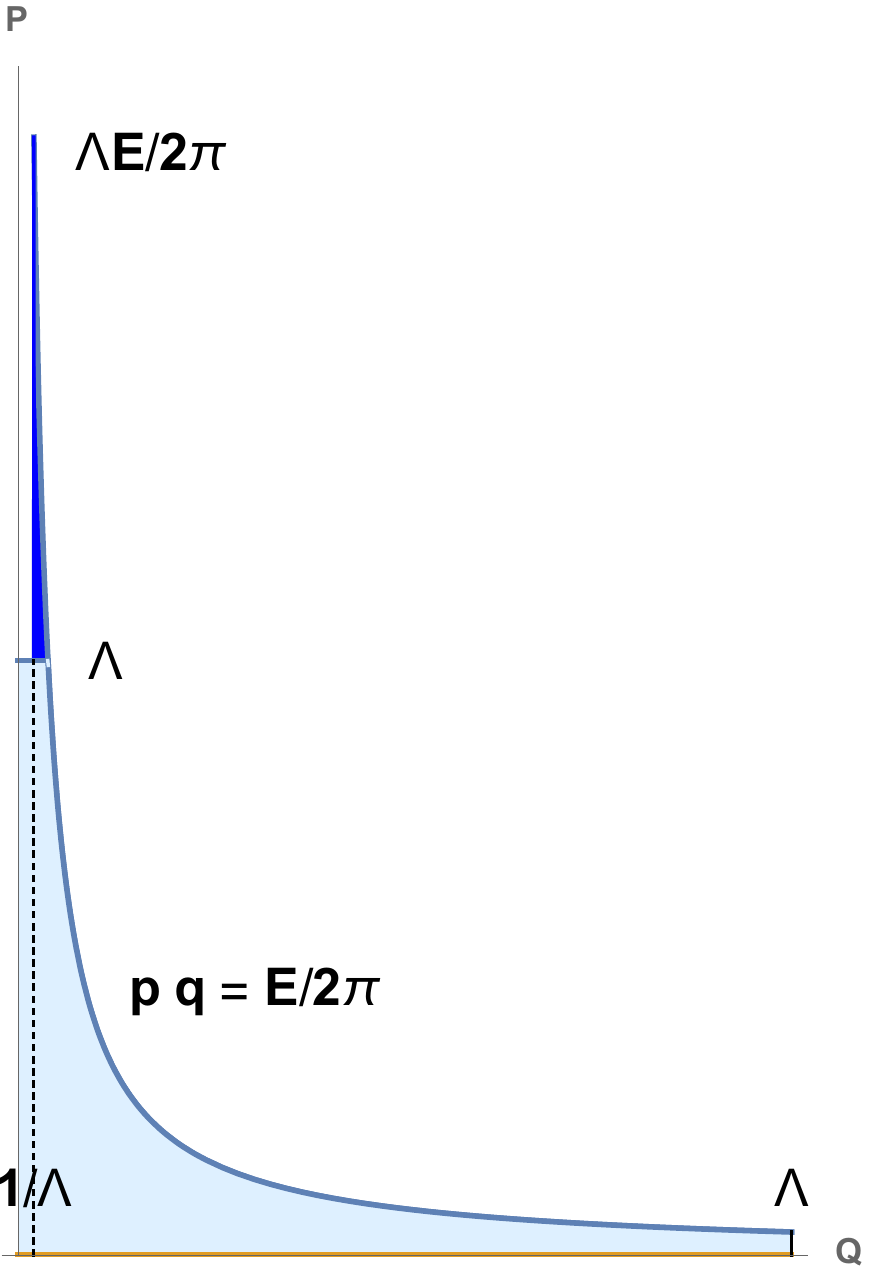}
\end{center}
\caption{\label{figfour} }
\end{figure}
Then there is also the rectangle, with vertices at $(0,0)$, $(1/\Lambda,0)$,$(1/\Lambda,\Lambda)$,$(0,\Lambda)$, in Figure \ref{fig2}, which does not fit inside Figure \ref{figthree}. It is  of area  equal to $1$ and thus corresponds to a single quantum state\footnote{This is reminiscent of the pole of the Riemann  zeta function, or equivalently to half of the boundary conditions
$
f(0)=0, \hat f(0)=0
$
that one needs to impose to use the Poisson Formula for the map $E$, as in \cite{Co-zeta}.}.

In order to compare this set-up with the Berry-Keating framework  we perform the symplectic transformation of the $(q,p)$ plane given by
$$
(q,p)\mapsto (\Lambda q,\frac{p}{\Lambda})=\sigma(q,p).
$$
This transformation preserves the area and the hyperbola $pq=\frac{E}{2\pi}$. It transforms Figure \ref{figfour} into the following Figure \ref{figfive}:

\begin{figure}[H]
\begin{center}
\includegraphics[scale=0.5]{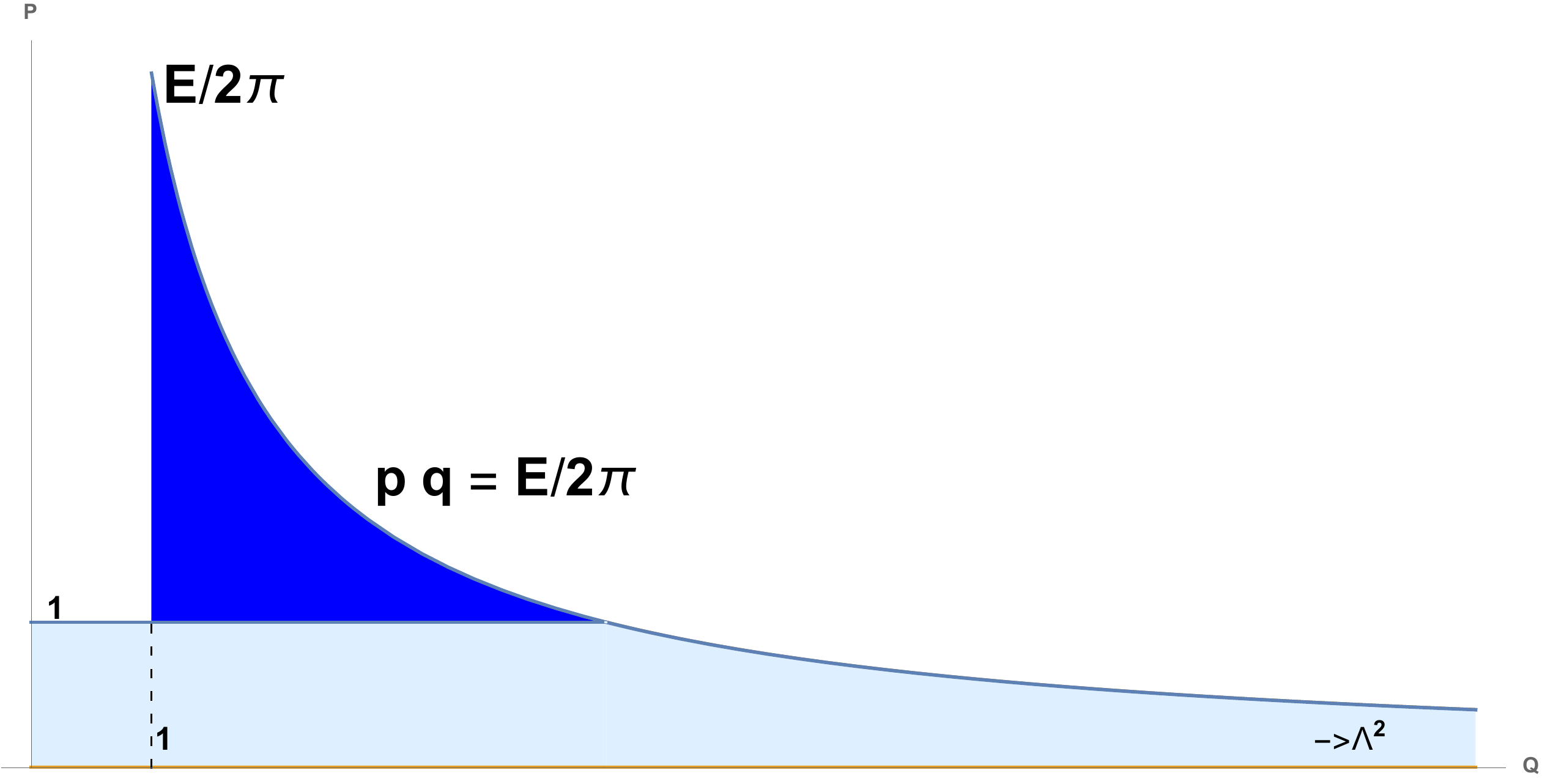}
\end{center}
\caption{\label{figfive} }
\end{figure}
One easily obtains that the missing (absorbed) part of the white light is given by the semiclassical picture of the ``Berry-Keating" Hamiltonian as in Figure \ref{figBK}.

This discussion shows that all the semiclassical computations for the Riemann  zeros viewed as an absorption spectrum (as in \cite{Co-zeta, CMbook}), are identical to those performed  with the ``Berry-Keating" Hamiltonian viewed as an emission spectrum (\cite{BK, BKe}). The presence of  the opposite signs involved in absorption and emission spectra does not create a conflict at this point. However, one should not conclude that the two models are equivalent.

\begin{figure}[H]
\begin{center}
\includegraphics[scale=0.38]{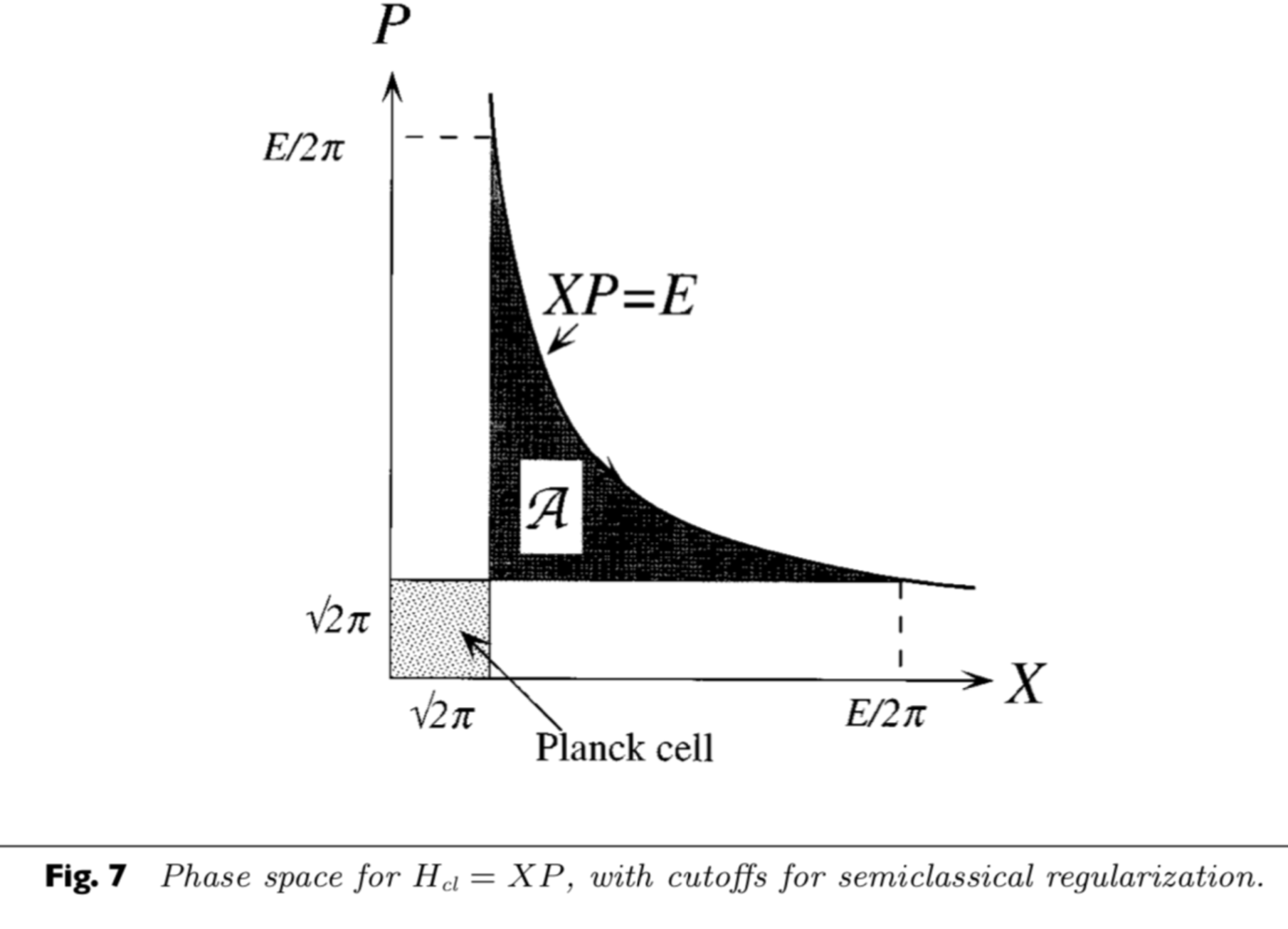}
\end{center}
\caption{Figure taken from \cite{BKe}, there is a different normalization of the Fourier transform which accounts for the discrepancy in the factors $2\pi$.\label{figBK} }
\end{figure}

 In fact, the minus sign coming from the absorption spectrum strikes back when one investigates the quantum fluctuations. In \cite{BK}, Berry and Keating are forced to set the Maslov phases all equal to $\pi$ in order to account for this minus sign. They write explicitly :
\begin{figure}[H]
\begin{center}
\includegraphics[scale=0.35]{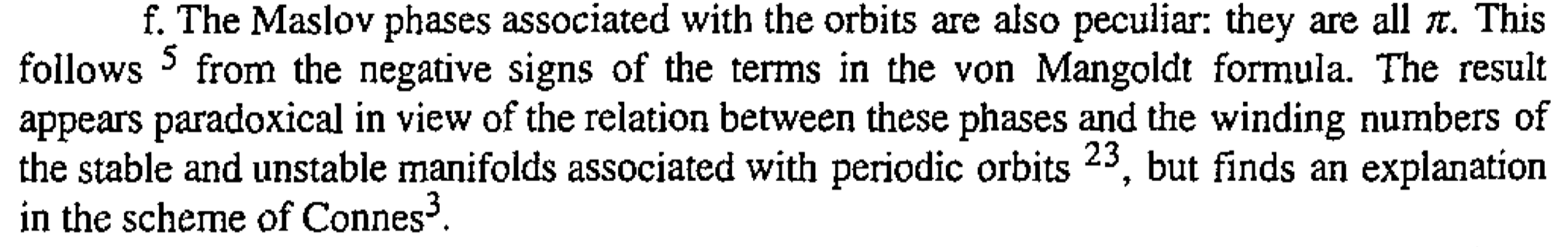}
\end{center}
\end{figure}

\section{SEMI-LOCAL TRACE FORMULA} \label{sectsemilocal} 

After the initial work \cite{Co-zeta} on the semi-local trace formula a better proof was obtained  based 
 on an idea of \cite{burnol} in the case of a single archimedean place. This development 
was explained in the class of the first author at Coll\`ege de France in 1999  \cite{Co99}.
We refer to Chapter 2 of \cite{CMbook} for a detailed exposition. A key tool here is the quantized calculus 
of Chapter IV of \cite{Co-book}.

\subsection{Quantized calculus}

 The main idea of the 
quantized calculus is to give an operator-theoretic version of the
calculus rules, based on the operator-theoretic differential
\begin{equation}\label{qdifferential}
\qd\, f := [F,f],
\end{equation}
where $f$ is an element in an involutive algebra $\cA$ represented
as bounded operators on a Hilbert space $\cH$. The  right-hand
side of \eqref{qdifferential} is the commutator with a self-adjoint
operator $F$ on $\cH$ with $F^2=1$. One then defines the analogue of differential 
forms of degree $n$ as linear combinations of monomials of the form
$$
\omega = f_0 \ \qd\, f_1\  \qd\, f_2 \ldots \qd\, f_n\in \Omega^n.
$$
The differential of an element of  $\Omega^n$ is defined as the graded variant of \eqref{qdifferential} 
\begin{equation}\label{qdifferential1}
\qd\, \omega := F \omega -(-1)^n \omega F.
\end{equation}
One checks, using $F^2=1$, that the square of the differential is $0$, $\ \qd^2 \omega=0$.
\medskip

Next, we briefly recall the framework for the quantized calculus in
one variable, as in Chapter  IV of \cite{Co-book}. We let functions $f$
of one real variable  act as multiplication operators in
$L^2(\R,ds)$, by
\begin{equation}\label{multoperf}
(f\,\xi)(s):= f(s)\,\xi(s)\,,\ \ \forall s\in \R, \ \xi \in L^2(\R,ds).
\end{equation}
We let $\fourier_{e_\R}$ be the Fourier transform with respect
to the  character  $e_\R(x):=e^{-2\pi ix}$ 
\begin{equation}\label{cFalpha}
\fourier_{e_\R}(\xi)(y):=\int_\R\,\xi(x)\,e^{-2\pi ixy}\,dy.
\end{equation}
We let ${\bf
1}_{[a,b]}$ be the characteristic function  of the interval $[a,b]\subset \R$ and  \begin{equation}\label{pr1}
\Pi_{[a,b]} := \fourier_{e_\R} \, {\bf 1}_{[a,b]}
\fourier_{e_\R}^{-1}
\end{equation}
be the conjugate by the Fourier transform $\fourier_{e_\R}$ of the
multiplication operator by ${\bf
1}_{[a,b]}$. 
Let  $H$ be the Hilbert transform $H = 2\; \Pi_{[0,\infty]} -1$
given, with care of the principal value, by
\begin{equation}\label{HilbTrans}
(H\,\xi)(s):= \frac{1}{i \pi} \int \frac{\xi(t)}{s-t} dt .
\end{equation}

\begin{definition}\label{qcalcdefn}
We define the {\em quantized differential} of $f$ to be the operator
\begin{equation}\label{qd}
\qd\, f:= [H,\,f]= H\,f-f\,H.
\end{equation}
\end{definition}

Thus, the quantized differential of $f$ is given by the kernel
\begin{equation}\label{qd1}
k(s,t)= \frac{i}{ \pi} \, \frac{f(s)-f(t)}{s-t}.
\end{equation}
 Definition
\ref{qcalcdefn} extends  to arbitrary modulated groups $C$, \ie to locally compact abelian groups endowed with a proper homomorphism $\Mod:C\to \R_+^*$. We let $\Mod(u)=\vert u\vert$ for $u \in C$. 
We let  $\widehat{C}$ be the
Pontrjagin dual of $C$ endowed with its Haar measure. The elements of $L^\infty(\widehat{C})$ act
as multiplication operators on the Hilbert space $\cH :=
L^2(\widehat{C})$. The operator $H$ on $\cH$ is given by
\begin{equation}\label{HilbTrsf}
 H:= 2 \fourier_C\,{\bf 1}_P\,\fourier_C^{-1} -1,
\end{equation}
where $\fourier_C: L^2(C) \to \cH$ is the Fourier transform  and 
${\bf 1}_P$ is the multiplication by the characteristic function of
the set $P=\{u \in C\,| \, \vert u\vert\geq 1\}$. Such $H$ plays the role of  the Hilbert transform.

In analogy with \eqref{qd}, we define the ``quantized" differential of
$f$ as the operator
\begin{equation}\label{qdd}
\qd\,f:= [H,\,f]= Hf-fH.
\end{equation}

We also use the following notation whose normalization fits with \eqref{pr1} :
\begin{equation}\label{PiabS}
\Pi_{[a,b]}:=  \fourier_C\,{\bf 1}_Y\,\fourier_C^{-1} \,,\ \ \
\text{ for } \ Y=\{u \in C\,| \,e^{2\pi a}\leq
 \vert u\vert\leq  e^{2\pi b}\}.
\end{equation}

\subsection{The Hilbert space $L^2
(X_S)$ for a finite set of places}
In the following we consider  the global field $\Q$ of rational numbers. In the ad\`ele ring of $\Q$ the additive Haar measure and the Haar measure of the multiplicative group of id\`eles are mutually singular. This problem does not arise when one restricts the attention to  a finite set $S$ of places of $\Q$.  We assume that $\infty \in S$. Let us consider
the locally compact ring
\begin{equation}\label{AQS}
 \A_{S}=\prod_{v\in S} \Q_v.
\end{equation}
It contains $\Q$ as a subring, using the diagonal embedding. We  let
$\Q_S$ denote the subring of $\Q$ given by all rational numbers whose
denominator only involves powers of primes $p \in S$. In other words,
\begin{equation}\label{AQSring}
 \Q_S=\{q\in \Q\,|\, |q|_v\leq 1\,,\ \forall v\notin S\}\,.
\end{equation}
The group $\Q^*_S$ of invertible elements of the ring $\Q_S$ is of
the form
\begin{equation}\label{GL1QS}
\Q_S^*=\GL_1(\Q_S)= \{ \pm p_1^{n_1} \cdots p_k^{n_k} \, :\,  p_j
\in S \setminus\{ \infty \} \,,\, n_j\in \Z\}.
\end{equation}

\medskip

We approximate the ad\`ele  class space by the semi-local ad\`ele  class space 
 $X_{S}$ which is the quotient
\begin{equation}\label{XQS}
X_{S}:=\A_{S}/\Q_S^* .
\end{equation}
The group
\begin{equation}\label{GL1AS}
C_{\Q,S}=\A_{S}^*/\Q_S^*, \  \ \ \A_{S}^*=\GL_1(\A_{S})= \prod_{p\in S} \GL_1(\Q_p)
\end{equation}
acts naturally by multiplication on the quotient $X_{S}$.

We  recall from \cite{Co-zeta,CMbook} that the Fourier transform $\fourier$ on the Bruhat-Schwartz space 
$\cS(\A_{S})$ induces a unitary operator on the Hilbert space
$L^2 (X_S)$. For each place $v\in S$ we choose a basic character
$\alpha_v$ of $\Q_v$. This makes it possible to identify the locally
compact abelian group $\Q_v$ with its Pontrjagin dual  by the pairing
\begin{equation}\label{PontrjaginPair}
\langle a, b\rangle := \alpha_v(a b), \ \ \  \forall a, b \in \Q_v.
\end{equation}
We normalize the additive Haar measure so that it is self-dual. Then
$\alpha =\prod_v \alpha_v$ is a  basic character of the additive
group $\A_S$, and we let $\fourier_\alpha$ denote the corresponding
Fourier transform, acting on the  space $\cS(\A_S)$. We refer to \cite{Co-zeta, CMbook} for the  proof of the following result
\begin{lemma} \label{fourierAQS}
Consider the character $\alpha =\prod_v \alpha_v$ and the
corresponding Fourier transform $ \fourier_\alpha   $ as above. The
map $f \mapsto
 \fourier_\alpha   (f)$, for $f \in \cS(\A_{S})$, extends uniquely to a unitary
operator $\underline\fourier_\alpha $ on the Hilbert space $L^2
(X_S)$.
\end{lemma}
The scaling operator $\urep (\lambda)$, for $\lambda \in \GL_1
(\A_{S})$, is defined as 
\begin{equation}\label{scalingS}
(\urep (\lambda) \, \xi) (a) = \xi (\lambda^{-1} a), \ \ \  \forall
a \in \A_{S}, \ \forall \xi\in \cS(\A_{S}).
\end{equation}
It is not unitary but its module is given by the module of the modulated group $C_{\Q,S}$
\begin{equation}\label{sc1}
\urep (\lambda)^* \urep (\lambda)= \urep (\lambda) \urep
(\lambda)^*= \vert \lambda \vert_S .
\end{equation}

We let $w$ be the unitary identification of $L^2(X_{S})$ with
$L^2(C_{\Q,S})$ (see \cite{CMbook} Proposition 2.30). It
intertwines the representation $\urep$  with the following 
regular representation $\vrep$ of $C_{\Q,S}$
\begin{equation}\label{vrepdefnfirstsemi}
(\vrep(\lambda)\,\xi)(v):=\,\xi(\lambda^{-1}\,v)\qqq \xi \in
L^2(C_{\Q,S})\,,
\end{equation}
 by the equality
\begin{equation}\label{wVhjsemi}
w \, \urep(\lambda)  w^{-1}= |\lambda|^{1/2}\,\vrep(\lambda) \qqq
\lambda\in C_{\Q,S}\,.
\end{equation}
For each
place $v$, the multiplicative group $\Q_v^*$ is a modulated group,
and by \cite{tate} there exists a unitary function
 $u_v \in L^\infty(\widehat{\Q_v^*})$ that allows one to express the Fourier transform $\fourier_{\alpha_v}$  relative to the
character $\alpha_v$, as a composition
of the inversion
$
I(f)(s):= f(s^{-1})
$ with the multiplicative convolution operator associated to $u_v$. The function $u_v$ is given by
the ratio of the local factors of $L$-functions (see \cite{tate}).
 This
classical formula   extends to our context as follows. We denote by $u$  the function
\begin{equation}\label{RieSiegelsemi}
u:= \prod_S\, u_v \circ \pi_v \,.
\end{equation}
where  $\pi_v:\widehat{C_{\Q,S}}\to \widehat{\Q_v^*}$ is the projection  dual to the embedding $\Q^*_v\hookrightarrow C_{\Q,S}$,  $x
\mapsto (1,1,\ldots ,x,\ldots ,1)$ which puts $x$ at the place $v$. We now have
\begin{lemma}(see \cite{CMbook}) \label{technfouriersemi}
On  $L^2(X_{S})$ one has
\begin{equation}\label{FwIPhisemi}
\underline\fourier_\alpha= w^{-1}\circ I \circ \fourier_C   ^{-1}
\circ u \circ   \fourier_C   \circ  w,
\end{equation}
where $\fourier_C$ is the Fourier transform of $C_{\Q,S}$ and  $
I(f)(s):= f(s^{-1})
$ is the inversion. 
\end{lemma}

\subsection{The semi-local trace formula}
As in the case of the single place $\infty$,  we introduce infrared
and ultraviolet cutoffs. For the infrared, we use the orthogonal
projection $P_\Lambda$ onto the subspace
\begin{equation}\label{subspPL}
\{ \xi \in L^2 (X_S) \, | \, \xi (x) = 0, \, \forall x \text{ with }
\vert x \vert > \Lambda \}.
\end{equation}
Thus, $P_\Lambda$ is the operator acting as multiplication by the
function $\rho_\Lambda$, with $\rho_\Lambda (x) = 1$ for $\vert x
\vert \leq \Lambda$ and $\rho_\Lambda (x) = 0$ for $\vert x \vert >
\Lambda$. 

We define an ultraviolet cutoff as $\widehat P_\Lambda =
\underline\fourier_\alpha P_\Lambda
\underline\fourier_\alpha^{-1}$, where $ \underline\fourier_\alpha
$ is the Fourier transform  of Lemma \ref{fourierAQS}, which depends
upon the choice of the basic character $\alpha =\prod_v \alpha_v$.
Then, we let
\begin{equation}\label{RLambdaS}
R_\Lambda = \widehat P_\Lambda  P_\Lambda.
\end{equation}
 We recall from Chapter 2 of \cite{CMbook} the
following result.

\begin{lemma}\label{techn1}
Let $S\ni \infty$ be as above. For any $\Lambda>0$ there is a unitary operator
\begin{equation}\label{WLambdaS}
W=W_\Lambda :L^2(X_{S})\to L^2(\widehat{C_{\Q,S}})
\end{equation}
such that, for any  $h_j\in \cS(C_{\Q,S})$, $j=1,2$, one has
\begin{equation}\label{WVh12}\begin{array}{l}
W\,\urep(\tilde h_1) \widehat P_\Lambda  P_\Lambda \urep(\tilde
h_2)W^*
 =\\[3mm]\qquad \qquad\hat h_1\left(\frac{1}{2} u^{-1}\,  \qd \,u\;
\Pi_{[-\infty,\frac{2\log \Lambda}{2\pi}]}+ \, \Pi_{[0,\frac{2\log
\Lambda}{2\pi}]}\right) \hat h_2.\end{array}
\end{equation}
Here $\tilde h_j(\lambda)=|\lambda|^{-1/2}h_j(\lambda)$, the
operator $\,\qd \,u$ is the quantized differential of the function
$u:= \prod_S\, u_v \circ \pi_v $ of \eqref{RieSiegelsemi} and $\hat
h_j$ is the multiplication operator by the Fourier transform
$\fourier_C(h_j)$.
\end{lemma}
This result together with the formula 
$$
(u_1u_2)^{-1}\,\qd (u_1u_2)=u_2^{-1}(u_1^{-1}\,\qd u_1) u_2+u_2^{-1}\,\qd u_2
$$ 
 allows one to compute $\Tr(\hat h_1 u^{-1}\,  \qd \,u\,\hat h_2) $ in \eqref{WVh12} as a sum, and get
\begin{equation}\label{semilocTrace0}
\Tr (\hat h_1\left(\frac{1}{2} u^{-1}\,  \qd \,u\right)\hat h_2) =  \sum_{v \in S}
\int'_{\Q^*_v} \frac{\vert w \vert^{1/2} }{ \vert 1-w \vert} \,h(w) d^* w,
\end{equation}
where $h=h_1 \star  h_2$ is such that $\vrep(h)=\vrep(  h_1)\vrep(h_2)$.
One then derives the following semi-local trace formula (see \cite{Co-zeta, CMbook})
\begin{theorem}\label{semilocalTr}
Let $\A_{S}$ be as above,   $\alpha=
\prod_{v\in S} \alpha_v$ a choice of basic character. Let $f \in \cS (C_{\Q,S})$ be a function
with compact support. Then, in the limit $\Lambda \to \infty$, one
has
\begin{equation}\label{semilocTrace}
\Tr (\urep(f)R_\Lambda) = 2f (1) \log \Lambda + \sum_{v \in S}
\int'_{\Q^*_v} \frac{f(w^{-1}) }{ \vert 1-w \vert} \, d^* w + o(1).
\end{equation}
\end{theorem}
Note that $\urep(f)=\int f(\lambda)\urep(\lambda)d^*\lambda$ corresponds by \eqref{wVhjsemi} to $\vrep(h)=\int h(\lambda)\vrep(\lambda)d^*\lambda$, where one has $h(\lambda)=|\lambda|^{1/2}f(\lambda)$. In terms of the test function $h$ the terms in the right hand side of \eqref{semilocTrace} take the more symmetric form 
\begin{equation}\label{semilocTracesym}
\int'_{\Q^*_v} \frac{f(w^{-1}) }{ \vert 1-w \vert} \, d^* w =\int'_{\Q^*_v} \frac{\vert w \vert^{1/2} }{ \vert 1-w \vert} \, h(w^{-1})\, d^* w =\int'_{\Q^*_v} \frac{\vert w \vert^{1/2} }{ \vert 1-w \vert} \, h(w)\, d^* w.
\end{equation}

\section{SONINE SPACES AND THE ATTEMPT OF X.-J.~LI}\label{sectsonin}

In  \cite{burnol1}, J. F. Burnol found a spectral realization of  the zeros of zeta, closely related  to the spectral realization of \cite{Co-zeta}, but using the Sonine spaces implemented  by L. de Brange  in his approach to RH. 
Motivated by the role of Sonine spaces,  as well as his own work on the semi-local trace formula \cite{XJL}, X.-J.~Li made, in early 2019, a  brave  attempt at proving the Weil positivity which is equivalent to RH (see \cite{EB}). Rather than  introducing a cutoff for large $q$ and $p$ as in \cite{Co-zeta}, he prescribes a cutoff near zero, \ie for small $q$ and $p$. In fact he does this in a balanced way \ie say requiring $q>\Lambda$ and $p>\Lambda^{-1}$. The coincidence between this prescription and the above  relation between the ``absorption" and "emission pictures" is very striking since the proposed cutoff near zero corresponds exactly to Figures \ref{figfour}, \ref{figfive}, \ref{figBK}. 

This attempt though fails. Nonetheless it is worthwhile to review it in some details and also to understand what is still missing, while retaining the interesting part. There are two basic facts on which Li's attempt is based: 
\begin{fact} \label{fact1} Let $A,B$ be positive bounded operators in a Hilbert space $\cH$ and assume that $AB$ is of trace class, then $\Tr(AB)\geq 0$.	
\end{fact}
\proof Let $A^{1/2}$ be the square root of $A$. Then, since in general the spectrum of a product $ST$ is the same as the spectrum of $TS$ except possibly for $0$, one obtains, denoting $\Spec^*$ the non-zero part of the spectrum
$$
\Spec^*(AB)= \Spec^*(A^{1/2}BA^{1/2})\geq 0.
$$
This shows that all eigenvalues of the trace class operator $AB$ are positive, and one concludes using Lidskii's theorem. \endproof 

\begin{fact} \label{fact2} The Riemann-Weil explicit formula applied to a test function $h \in \cS (C_{\Q})$ with compact support only involves finitely many places.	
\end{fact}
\proof By hypothesis the function $h$ only depends on the modulus. The contribution of each prime $p$ only involves the value of  $h$ on non-zero powers of $p$ so that only finitely many primes contribute. \endproof 
The natural idea then, is to use the semi-local framework of \cite{Co-zeta} recalled in Section 
\ref{sectsemilocal} and more precisely  formula \eqref{semilocTrace0} to show an inequality of the form 
\begin{equation}\label{semilocTracex}
  \sum_{v \in S}
\int'_{\Q^*_v} \frac{\vert w \vert^{1/2} }{ \vert 1-w \vert} \,h(w) d^* u \leq 0, \ \ \  \  h=h_1 \star  h_2, \ h_2(w):=\overline {h_1(w^{-1}}).
\end{equation} 
Note that by construction the operators $\hat
h_j$ given as the multiplication operators by the Fourier transform
$\fourier_C(h_j)$ fulfill, when $h_2=h_1^*$ as in \eqref{semilocTracex}, the following 
\begin{fact} \label{fact3} When $h_2=h_1^*$, the operator $\hat h_1 \hat h_2$ is positive.
\end{fact}

Now, using the cyclic property of the trace and the commutativity of the algebra of multiplication by functions, one gets from  \eqref{semilocTrace0} the following equality
\begin{equation}\label{semilocTracey}
  \sum_{v \in S}
\int'_{\Q^*_v} \frac{\vert w \vert^{1/2} }{ \vert 1-w \vert} \,h(w) d^* w = \Tr \left(\hat h_1 \hat h_2(\frac{1}{2} u^{-1}\,  \qd \,u)\right).
\end{equation}
The idea of X.-J.~Li is to replace the operator $R_\Lambda=\widehat P_\Lambda  P_\Lambda$ by the  positive operator
$
T=P\hat{ P}P,  
$ where $P$ is the projection corresponding to the cutoff near zero, \ie restricting to $\vert w \vert \geq 1$. One gets, using the notations of \eqref{semilocTrace0},   that 
$$
T=Pu^*(1-P)uP.
$$
One sees that Fact \ref{fact1} would imply \eqref{semilocTracex} provided one could prove that  $$Pu^*(1-P)uP=-\frac 12 	u^{-1}\,  \qd \,u.$$ As we shall see this is equivalent to  controlling the sign of the operator $u^{-1}\,  \qd \,u$. The quantized calculus gives in fact the following general criterion to control the sign of the ``logarithmic derivative" 
$u^{-1}\,  \qd \,u$ of a unitary $u$.
\begin{lemma} \label{semiloclem} Let $\cH$ be a Hilbert space,  $u$  a unitary operator and $F=2P-1$, where $P$ is an orthogonal projection. Then, with $\, \qd u:=[F,u]$,   the following three conditions are equivalent
\begin{enumerate}
\item $Pu^*(1-P)uP=-\frac 12 	u^{-1}\,  \qd \,u$
\item $u^{-1}\,  \qd \,u \leq 0$
\item $Pu=PuP$
\end{enumerate}
\end{lemma}	
\proof $(i)$ $\Rightarrow$ $(ii)$ since $Pu^*(1-P)uP$ is a positive operator. 

$(ii)$ $\Rightarrow$ $(iii)$. If $u^{-1}\,  \qd \,u \leq 0$ one has $u^*[P,u]\leq 0$ and $u^*Pu \leq P$, thus $Pu^*P=u^*P$  since the range of $u^*P$ is inside $P$. Thus, taking adjoints, one obtains $Pu=PuP$. 

\begin{figure}[H]	\begin{center}
\includegraphics[scale=0.8]{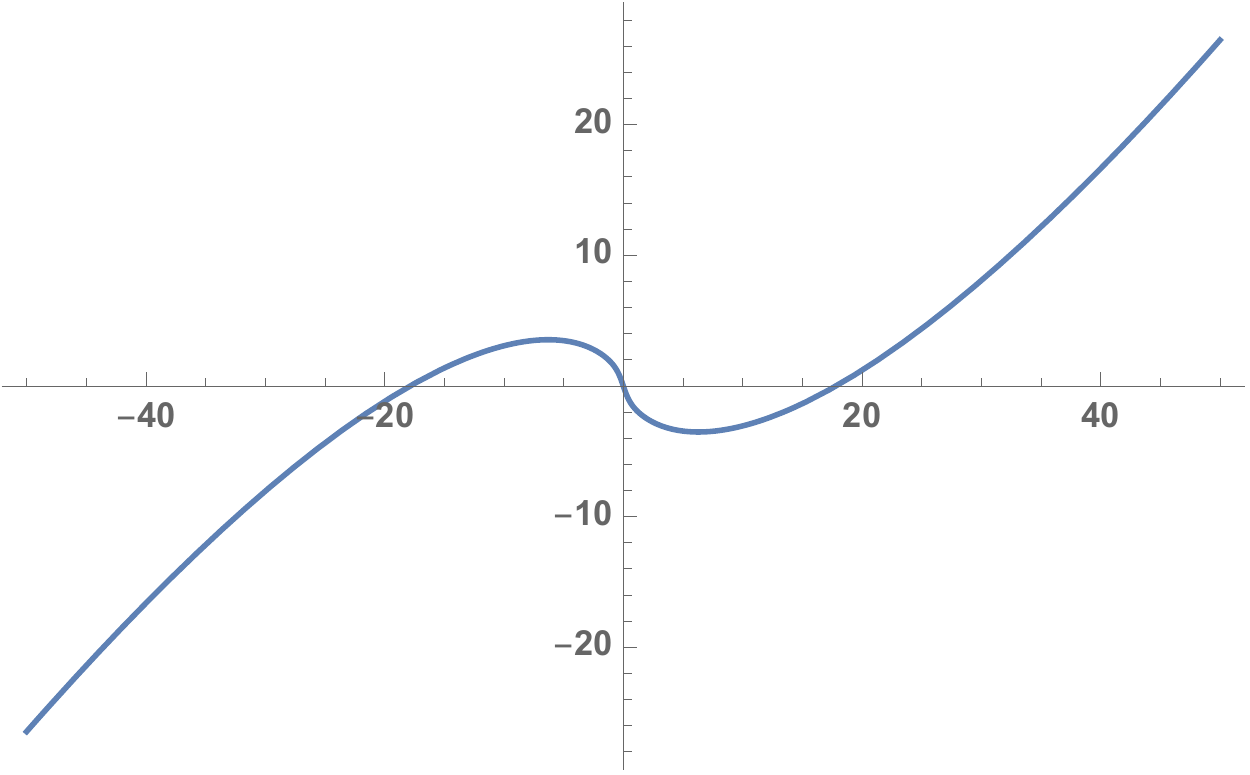}
\end{center}
\caption{Riemann-Siegel function. \label{riemannsiegel} }
\end{figure}

$(iii)$ $\Rightarrow$ $(i)$. Assume $Pu=PuP$. Taking adjoints one derives $Pu^*P=u^*P$, and thus the projection $Q=u^*Pu$ is less than $P$, one has $PQP=Q$ and  
$$
Pu^*(1-P)uP=P(1-Q)P=P-Q=-u^*[P,u]=-\frac 12 	u^{-1}\,  \qd \,u
$$
which proves $(i)$.
\endproof
When working with the quantized calculus in one variable one gets the following
\begin{corollary} \label{semiloccor} Let $\cH$ and $F$ be as in Definition \ref{qcalcdefn}, and $u$ a unitary function. Then the equivalent conditions of Lemma \ref{semiloclem} hold  if and only if $u^*$ is an inner function in the sense of Beurling \cite{Beu}. 	
\end{corollary}
\proof The function $u^*$ is inner if and only if it preserves the subspace which is the range of the projection $P$, with $F=2P-1$. Then Lemma \ref{semiloclem} gives the required equivalence. \endproof 
Equipped with  Lemma \ref{semiloclem} and Corollary \ref{semiloccor}, we can now understand why the direct approach proposed by X.-J.~ Li fails to prove Weil's positivity. Indeed, in our case and taking a single place to simplify the discussion, the issue is whether the functions $u_v^*$, where $u_v$  is given by the ratio of local factors on the critical line, happen to be  inner functions. In fact we can ask if $u_v$ 
 belongs to one of the 
Hardy spaces $\H^\infty(\C_{\frac 12}^\pm)$ of the half plane $$\C_{\frac 12}^\pm:=\{z\in \C\mid \pm\Re(z)>1/2\}.$$ 
 At a finite prime $p$ the function $u_p$ is given by the ratio of local factors 
\begin{equation}
u_p(s):= \frac{1-p^{-(1-z)}}{1-p^{-z}}, \ \ z=\frac 12 +is. 	
\end{equation}
It is a function of modulus $1$ for $s\in \R$.  The zeros of the denominator $1-p^{-z}$ give poles,
$$
1-p^{-z}=0\iff z\in \frac{2\pi i}{\log p}\Z
$$
which show that in the left half-plane the function is not holomorphic. On the other hand in the right half-plane the denominator does not vanish and the function is holomorphic. However it is not bounded since for real values of $z$ it behaves as follows
$$
\frac{1-p^{-(1-z)}}{1-p^{-z}}\sim -p^{z-1}, \ \ z\to \infty. 
$$
We now look at the archimedean place. There, the ratio of local factors is 
\begin{equation}\label{archimloc}
u_\infty(s):=\phi(\frac 12 +is), \ \ \phi(z)=\frac{\pi ^{\frac{1-z}{2}-\frac{z}{2}} \Gamma \left(\frac{z}{2}\right)}{\Gamma \left(\frac{1-z}{2}\right)}.
\end{equation}
The function $\phi(z)$ is of modulus $1$ for $z=\frac 12+is$, and has poles in the left half-plane.  It is holomorphic in the right half-plane $\Re(z)\geq \frac 12$, however it is not bounded there. In fact one has using the complement formula
$$
\Gamma(x)\Gamma(1-x)=\frac{\pi}{\sin(\pi x)}\Rightarrow 
\frac{\pi}{\Gamma(1-x)}=\Gamma(x)\sin(\pi x)
$$
with $x=\frac 12(1+z)$, that
$$
\phi(z)=\pi ^{-\frac{1}{2}-z} \Gamma \left(\frac{z}{2}\right)\Gamma \left(\frac{1+z}{2}\right)\sin\left(\pi \frac{1+z}{2}\right).
$$
Thus even on the real axis the function $\phi(z)$ that vanishes at odd integers takes very large values on even integers such as 
$$
\phi(20)=\frac{1856156927625}{8 \pi ^{20}}\sim 26.4562.
$$
Thus we conclude that none of the functions $u_v(s)^*$ fulfills the requirement of Corollary \ref{semiloccor}.
 We finally remark that, a priori,  the inequality \eqref{semilocTracex} could still hold even though the above described attempt to prove it fails. In  fact one has the following
\begin{fact} \label{fact4} The inequality \eqref{semilocTracex} does not hold in general.
\end{fact}
\proof We consider 
 the simplest case when S is reduced to the single archimedean place.
Inequality \eqref{semilocTracex} would then imply, using the expression of the left hand side in terms of the logarithmic derivative of $u_\infty(s)$,  that  this logarithmic derivative (multiplied by $i$)
 has a constant sign. But $u_\infty(s)=e^{2i\,\theta(s)}$ where $\theta(s)$ is the Riemann-Siegel angular function (see \cite{CMbook} Chapter II, \S 5.1, Lemma 2.20). Thus the inequality \eqref{semilocTracex}   would imply 
that the Riemann-Siegel angular function $\theta(s)$ is a monotonic function, but
by looking 
at its graph  (Figure  \ref{riemannsiegel}) one sees that it is not the case.

\section{AN OPERATOR THEORETIC PROBLEM} \label{sectoptheory}

By a result of A. Weil (see  \cite{EB})  RH is equivalent to  the negativity of the right hand side of the explicit formulas for all test  functions $f$ of the form 
$$
f(x)=\int_0^\infty g(xy)\overline {g(y)}dy, \ \int_0^\infty g(x)dx=0, \ \int_0^\infty g(x)d^*x=0.
$$
One may moreover restrict attention to test functions $g$ with compact support (see \cite{B3}). With these notations,  the right hand side of the explicit formulas  corresponds, in the more general case of $L$-functions,  to the terms in \eqref{semilocTrace}. To pass to the more symmetric form \eqref{semilocTracesym}, one then lets $h_1(x)=x^{1/2}g(x)$,  thus getting 
$$
f(x)=\int_0^\infty g(xy)\overline {g(y)}dy=x^{-1/2}\int_0^\infty h_1(xy)\overline {h_1(y)}d^*y.
$$
Thus, with $h_2=h_1^*$, so that $h_2(x):=\overline{h_1(x^{-1})}$ one obtains 
$$
h(x)=x^{1/2}f(x)=\int_0^\infty h_1(xy)\overline {h_1(y)}d^*y
=\int_0^\infty h_1(xy^{-1})\overline {h_1(y^{-1})}d^*y= (h_1\star h_2)(x).
$$
Hence  Weil's positivity corresponds well to \eqref{semilocTracex} but one needs to impose the two conditions 
\begin{equation}\label{h0h2}
	\int_0^\infty h_1(x)x^{\pm 1/2}d^*x=0.
	\end{equation}
The above discussion suggests the following
\begin{guess} \label{conjsl} The semi-local operator theoretic framework with $S:=\{\infty\}\cup \{p \mid p< q\}$ suffices to prove the Weil inequality for all test functions with support in the interval $(q^{-1/2},q^{1/2})$.	
\end{guess}

\subsection{Poisson Formula and basic characters}

The Poisson summation formula gives the conceptual reason why the ratio of local factors on the critical line appear when writing the Fourier transform as the composition of inversion with a multiplicative convolution. To understand this point we start with a single place, \ie we take $S=\{\infty\}$, and proceed at the formal level not caring too much about some technical points. One introduces the following operator acting on even functions $f\in \cS(\R)$ satisfying $f(0)=\hat f (0)=0$,  
\begin{equation}\label{mapE}
	\mapE(f)(x):=\vert x\vert^{1/2}\sum_{n>0} f(nx).
		\end{equation}
The Poisson summation formula then asserts that, provided one takes the Fourier transform $\fourier_{e_\R}$ associated to the basic character $e_\R$ for which the lattice $\Z\subset \R$ is self-dual (\ie the basic character $e_\R(x)=e^{-2\pi ix}$
 used in \eqref{cFalpha}), one has 
 \begin{equation}\label{Poisson}
	\mapE(\fourier_{e_\R}(f))(x)= \mapE(f)(x^{-1})\qqq x\in \R_+^*.
		\end{equation}
		We define  the duality $\langle \R^*_+, \R \rangle$ by the
bicharacter
\begin{equation}\label{FwIPhi0}
\mu(v,s)=v^{-is} \,,\ \ \forall v\in \R^*_+, s\in \R\,,
\end{equation}
 so that the Fourier transform\index{Fourier transform} $ \fourier_\mu   :
L^2(\R^*_+)\to L^2(\R)$  is given by
\begin{equation}\label{PhiFourier}
 \fourier_\mu   (f)(s):=\int_0^{\infty} f(v) v^{-is}d^*v\,.
\end{equation}
	When reading \eqref{mapE} in the dual of the multiplicative group $\R_+^*$, the multiplication of the variable $x\mapsto nx$ which is a translation in $\R_+^*$ becomes the multiplication of the function by $n^{is}$.  With $f_n(x):= x^{1/2}f(nx)$ one then  has
	$$
	\fourier_\mu(f_n)(s)= \int_0^{\infty} f(nv) v^{1/2-is}d^*v=n^{-(1/2-is)}\fourier_\mu(f)(s),  
	$$
so that, at this formal level, we obtain 
$$
\fourier_\mu(\mapE(f))(s)=\zeta(1/2-is)\fourier_\mu(f)(s).
$$
Thus we may think of $\mapE$, when read in Fourier, as multiplication by $\zeta(1/2-is)$. The fact that $\mapE$ can be viewed, in Fourier, as a multiplication operator is not surprising since by construction $\mapE$ commutes with scaling. Another map which commutes with scaling is the composite $I \circ \fourier_{e_\R}$ of the inversion $I$ with the additive Fourier transform $\fourier_{e_\R}$. Thus this composition is given,  when read in Fourier, as the multiplication by a function $u$ which is  of modulus $1$ due to unitarity.  The point then is that the Poisson Formula \eqref{Poisson} states that the map $\mapE$ conjugates  $\fourier_{e_\R}$ with the inversion $I$ and one obtains, still at this formal level, that 
$$
u(s)=\zeta(1/2-is)/\zeta(1/2+is).
$$
By the functional equation this quotient is the ratio \eqref{archimloc} of archimedean local factors. 

This discussion adapts to the semi-local case provided one replaces the summation over $\Z$ in the Poisson Formula by the summation over the discrete subgroup \eqref{AQSring}. 
Since $\Q_S$ is a ring, this summation is invariant under multiplication by the group $\Q_S^*$ of units. Moreover, every element $q\in \Q_S$ is uniquely a product $q=u m$, where $u \in\Q_S^*$ is a unit and $m\in M_S$ is an element of the multiplicative monoid of positive integers prime to all $p\in S$. The role of the map $\mapE$ of \eqref{mapE} is now played  by
\begin{equation}\label{mapEsl}
	\mapE(f)(x):=\vert x\vert^{1/2}\sum_{m\in M_S} f(mx).
		\end{equation}
		When computed in Fourier as above but now involving characters of the group $C_{\Q,S}$, \eqref{mapEsl} produces sums for $L$ functions which reduce, for the trivial character, to 
		$$
		\sum_{M_S} m^{-1/2+is}=\zeta(1/2-is)\ \prod_{p\in S}(1-p^{-1/2+is}).
		$$
	Thus  we see, heuristically as above, that the Poisson Formula gives the conceptual explanation for Lemma \ref{technfouriersemi}. 
There is in general a choice involved for basic additive characters $\alpha_S$ of $\A_S$ such that $\Q_S$ is a self-dual lattice, but this ambiguity disappears when working in the quotient $X_{S}=\A_{S}/\Q_S^* $. In particular, the Fourier transform $\underline\fourier_\alpha$ is independent of such choices. Indeed, the Fourier transforms of the same function for different choices of characters for which  $\Q_S$ is a self-dual lattice are related by the action (by scaling) of the group $\Q_S^*$ and thus their difference vanishes in the space $L^2 (X_S)$. The Weil inequality needs the above normalization of additive characters and this indicates that the Poisson Formula should play a key role in the solution of the Conjecture \ref{conjsl}. 

In fact, there is another striking reason why the Poisson Formula is intimately related to  
Conjecture \ref{conjsl}. It is the construction, due to J.P. Kahane (see \cite{burnol2}), of non-trivial elements of the Sonine space, \ie of functions which vanish identically near the origin as well as their Fourier transform. Let $\Pi:=\sum_{n\in \Z}\delta_n$ be the Poisson distribution, which is the sum of Dirac distributions on the integers. Then the tempered distribution $\cD(\Pi)$ has the above vanishing property, when $\cD$ is a  differential operator  whose range $\cD(\cS(\R))$ consists of functions which fulfill the two conditions  
\begin{equation}\label{twocond}
	f(0)=\hat f(0)=0.
		\end{equation}
		The simplest choice for $\cD$ is $\cD(f)(x):=xf''(x)$, since one checks by integration by parts that the integral of $\cD(f)$ vanishes. For any such choice of $\cD$ the distribution $\cD(\Pi)$ has the required Sonine support property. One then obtains a solution in Schwartz space by applying to $\cD(\Pi)$  the smoothing by convolution with a smooth function $\phi$ with small support near zero, and multiplication by the Fourier transform of $\phi$. In fact the above choice of the differential operator $\cD$, though simple, is not symmetric with respect to Fourier transform.  One checks that the simple function $\cD:=D_u^2+D_u$ of our scaling Hamiltonian $D_u(f)x):=xf'(x)$ does indeed work, commutes with Fourier transform, and one gets 
		 \begin{equation}\label{kahane}
	\cD(f)(x)= x^2 f''(x)+2 xf'(x).
		\end{equation}
One checks directly that the range $\cD(\cS(\R))$ consists of functions which fulfill the two conditions \eqref{twocond}. All this construction extends verbatim to the semi-local case.

\begin{remark}\label{tropical}  
	In an intriguing manner the adjoint expression $\Delta_2:=D_u^2-D_u$ provides the tropical structure of the Scaling Site as explained in \cite{CCscal1}, equation (1).
\end{remark}
 
\subsection{Time versus Energy}
In the approach to RH advocated by M. Berry and J. Keating \cite{BK, BKe}, one looks for an interpretation of the zeros of the Riemann zeta function as energy levels of a quantum system. In our own approach however the zeros appear naturally as ``times" \ie in a dual manner (see \cite{Co-zeta}). While this fact might look perplexing at first sight, we shall briefly explain here below why this dual point of view is in fact more natural in view of the generalization of the Riemann zeta function in the framework of global fields. Indeed, in this set-up one meets  much simpler avatars of the Riemann zeta function. They are associated to a curve $C$ over a finite field $\F_q$. It turns out that these analogues of the Riemann zeta function $\zeta(s)$ are in fact functions of the form $L(q^{-s})$, where $q$ is the cardinality of the finite field over which the curve is defined. Moreover, $L(z)$ is a rational fraction and its zeros  are determined by the zeros $z_j$ of the polynomial numerator $P(z)$ whose degree is twice the genus $g$ of the curve. By a famous theorem of A. Weil, all  $z_j$  are on the circle of radius $q^{-1/2}$. Thus the zeros $s$ of $L(q^{-s})$ are all of real part $\Re(s)=\frac 12$ and their imaginary parts $\Im(s)$ are of the form 
$$
\Im(s)=\frac{1}{\log q}\left(\alpha_j+2\pi k \right), \ j\in \{1, \ldots 2g\}, \ k\in \Z,
$$
where the $-\alpha_j$ are the arguments of the $z_j$ and are determined only modulo $2\pi$. This type of distribution of numbers, being periodic,  is very natural as a distribution in ``time", not in ``energy".
\subsection{Analysis and Geometry}

The above discussion does not dismiss the nagging thought we started with. In fact there is a real possibility that one could find a way to solve Conjecture \ref{conjsl}, by first concentrating on the simplest case of a single place and finding a way to make use of the support condition on the test function. A scary aspect of this strategy is that  Weil's positivity is very much dependent on the careful choices of principal values and it is difficult to avoid computational mistakes. This analytical approach is in this respect, by its vulnerability, quite different from the geometric approach that we have gradually developed in \cite{Crh,CCsurvey}. While the ad\`ele class space is clearly implicit in the semi-local analytic approach described above, its geometric structure is only unveiled at the global level and it reveals its conceptual meaning as the points of the  Scaling Site, where the action of the scaling group reveals itself as the action of the Frobenius.







\end{document}